\documentclass[aos,citesort,seceqn,dvips]{arximspdf}
\usepackage{graphicx}

\doi{10.1214/10-AOS859}
\volume{39}
\issue{2}
\pubyear{2011}
\firstpage{865}
\lastpage{886}

\makeatletter
\newcommand{\pa}{\operatorname{pa}} % for parents
\newcommand{\rank}{\operatorname{rank}}
\newcommand{\kernel}{\operatorname{kernel}}
\newproclaim{definition}{Definition}
\newproclaim{remark}{Remark}
\newproclaim{question}{Question}
\newtheorem{theorem}{Theorem}
\newproclaim{example}{Example}
\newtheorem{corollary}{Corollary}
\newtheorem{lemma}{Lemma}
\newtheorem{proposition}{Proposition}
\makeatother

\begin{document}
\begin{frontmatter}

\title{Global identifiability of linear structural equation models}
\runtitle{Global identifiability of linear structural equation models}

\begin{aug}
\author[A]{\fnms{Mathias} \snm{Drton}\corref{}\ead[label=e1]{drton@uchicago.edu}\thanksref{a1}},
\author[A]{\fnms{Rina} \snm{Foygel}\ead[label=e2]{rina@uchicago.edu}}
\and
\author[B]{\fnms{Seth} \snm{Sullivant}\ead[label=e3]{smsulli2@ncsu.edu}\thanksref{a2}}
\runauthor{M. Drton, R. Foygel and S. Sullivant}
\affiliation{University of Chicago,  University of Chicago\\
and~North~Carolina~State~University}
\address[A]{M. Drton\\
 R. Foygel\\
Department of Statistics\\
 University of Chicago\\
Chicago, Illinois\\
USA\\
\printead{e1}\\
\hphantom{\textsc{E-mail}: }\printead*{e2}} %adresu isvedimo komanda
%gale!
\address[B]{S. Sullivant\\
Department of Mathematics\\
North Carolina State University\\
Raleigh, North Carolina\\
USA\\
\printead{e3}}
\end{aug}
\thankstext{a1}{Supported by NSF Grant DMS-07-46265 and an Alfred P. Sloan
Fellowship.}
\thankstext{a2}{Supported by NSF Grant
DMS-08-40795 and the David and Lucille Packard Foundation.}

% HISTORY:
\received{\smonth{3} \syear{2010}}
\revised{\smonth{9} \syear{2010}}

% ABSTRACT
%
\begin{abstract}
Structural equation models are multivariate statistical models that are
defined by specifying noisy functional relationships among random
variables. We consider the classical case of linear relationships and
additive Gaussian noise terms. We give a necessary and sufficient
condition for global identifiability of the model in terms of a mixed
graph encoding the linear structural equations and the correlation
structure of the error terms. Global identifiability is understood to
mean injectivity of the parametrization of the model and is fundamental
in particular for applicability of standard statistical \mbox{methodology}.
\end{abstract}

% KEYWORDS
%

\end{frontmatter}

%s1 ###
\section{Introduction}
\label{sec:intro}

A mixed graph is a triple $G=(V,D,B)$ where $V$ is a finite set of nodes
and $D,B\subseteq V\times V$ are two sets of edges. The edges in $D$ are
directed, that is, $(i,j)\in D$ does not imply $(j,i)\in D$. We denote and
draw such an edge as $i\to j$. The edges in $B$ have no orientation; they
satisfy $(i,j)\in B$ if and only if $(j,i)\in B$. Following tradition in
the field, we refer to these edges as bidirected and denote and draw them
as $i\leftrightarrow j$. (In figures, we will draw bidirected edges
also as dashed
edges for better visual distinction.) We emphasize that in this setup the
bidirected part $(V,B)$ is always a simple graph, that is, at most one
bidirected edge may join a pair of nodes. Moreover, neither the bidirected
part $(V,B)$ or the directed $(V,D)$ contain self-loops, that is,
$(i,i)\notin D\cup B$ for all $i\in V$. In the main part of this work,
the considered mixed graphs are acyclic, which means that the directed part
$(V,D)$ is a directed graph without directed cycles.

Enumerate the vertex set as $V=[m]:=\{1,\ldots ,m\}$. Let $\mathbb
{R}^D$ be
the set of matrices $\Lambda=(\lambda_{ij})\in\mathbb{R}^{m\times
m}$ with
$\lambda_{ij}=0$ if $i\to j$ is not in\vspace*{-1pt} $D$. Write
$\mathbb{R}^D_{\mathrm{reg}}$ for the subset of matrices
$\Lambda\in\mathbb{R}^D$ for which $I-\Lambda$ is invertible, where $I$
denotes the identity matrix. Let $\mathit{PD}(m)$ be the cone of positive
definite $m\times m$ matrices. Define $\mathit{PD}(B)$ to be the set of
matrices $\Omega=(\omega_{ij})\in\mathit{PD}(m)$ with $\omega
_{ij}=0$ if
$i\not= j$ and $i\leftrightarrow j$ is not an edge in $B$. Write
$\mathcal{N}_m(\mu,\Sigma)$ for the multivariate normal distribution with
mean $\mu\in\mathbb{R}^m$ and covariance matrix~$\Sigma$.

\begin{definition}
\label{def:sem}
The linear structural equation model $\mathcal{M}(G)$ associated
with an acyclic mixed graph $G=(V,D,B)$ is the family of
multivariate normal distributions $\mathcal{N}_m(0,\Sigma)$ with
\[
\Sigma= (I-\Lambda)^{-T}\Omega(I-\Lambda)^{-1}
\]
for $\Lambda\in\mathbb{R}^D_{\mathrm{reg}}$ and
$\Omega\in\mathit{PD}(B)$.
\end{definition}

The set of parents of a node $i$, denoted $\pa(i)$, comprises the
nodes $j$ with $j\to i$ in $D$. The graphical model just defined
is most naturally motivated in terms of a system of linear structural
equations:
%
%e1.1 ###
\begin{equation}
\label{eq:sem}
Y_j = \sum_{i\in\pa(j)} \lambda_{ij} Y_i + \varepsilon_j,
\qquad j=1,\ldots ,m.
\end{equation}
If $\varepsilon=(\varepsilon_1,\ldots ,\varepsilon_m)$ is a random
vector following the multivariate normal distribution
$\mathcal{N}(0,\Omega)$ and $\Lambda\in\mathbb{R}^D_{\mathrm{reg}}$, then
the random vector $Y=(Y_1,\ldots ,Y_m)$ is well defined as a solution to
the equation system in (\ref{eq:sem}) and follows a centered
multivariate normal distribution with covariance matrix
$(I-\Lambda)^{-T}\Omega(I-\Lambda)^{-1}$.

\begin{remark}
Assuming centered distributions presents no loss of generality. An
arbitrary mean vector could be incorporated by adding an intercept
constant $\lambda_{i0}$ to each equation in (\ref{eq:sem}). The results
discussed below would apply unchanged.
\end{remark}

Linear structural equation models are ubiquitous in many applied fields,
most notably in the social sciences where the models have a long tradition.
Recent renewed interest in the models stems from their causal
interpretability; compare \cite{spirtes2000,pearl2009}. While current
research is often concerned with non-Gaussian generalizations of the
models, there remain important open problems about the linear Gaussian
models from Definition~\ref{def:sem}. These include the following
fundamental problem, which concerns the global identifiability of the model
parameters.

\begin{question}
\label{que:ident}
For which mixed graphs $G=(V,D,B)$ is the rational~para\-metrization
\[
\phi_G\dvtx  (\Lambda,\Omega) \mapsto(I-\Lambda)^{-T}\Omega(I-\Lambda)^{-1}
\]
an injective map from $\mathbb{R}^D_{\mathrm{reg}}\times\mathit{PD}(B)$
to the
positive definite cone $\mathit{PD}(m)$?
\end{question}

According to our first theorem, proven later on in
Section~\ref{sec:cyclic-models}, we can restrict attention to acyclic
mixed graphs.

\begin{theorem}
\label{thm:acyclic}
If $G$ is a mixed graph for which the parametrization $\phi_G$ is
injective, then $G$ is acyclic.
\end{theorem}

The nodes of an acyclic mixed graph $G=(V,D,B)$ can be ordered
topologically such that $i\to j\in D$ only if $i<j$. Under a topological
ordering of the nodes, all matrices in $\mathbb{R}^D$ are strictly
upper-triangular. Hence, $\mathbb{R}^D_{\mathrm{reg}}=\mathbb{R}^D$ because
$\det(I-\Lambda)=1$ for all $\Lambda\in\mathbb{R}^D$. Moreover, the
parametrization $\phi_G$ is a polynomial map in the entries of
$\Lambda$
and $\Omega$ when $G$ is acyclic.

Characterizing the graphs with injective parametrization is important
because failure of injectivity can lead to failure of standard statistical
methods. We briefly exemplify this issue for the models considered here
and point the reader to \cite{drtonlrt} and references therein for a more
detailed discussion. Briefly put, the problem is due to the fact that
failure of injectivity can result in parameter spaces that are not smooth
manifolds; compare in particular the examples in Section~1 of
\cite{drtonlrt}.

%
%f1 ###
\begin{figure}

\includegraphics{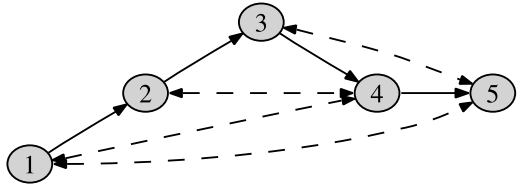}

\caption{Acyclic mixed graph inducing a singular model.}
\label{fig:graph-sim}
\end{figure}

%
%f2 ###
\begin{figure}
\centering
\begin{tabular}{cc}

\includegraphics{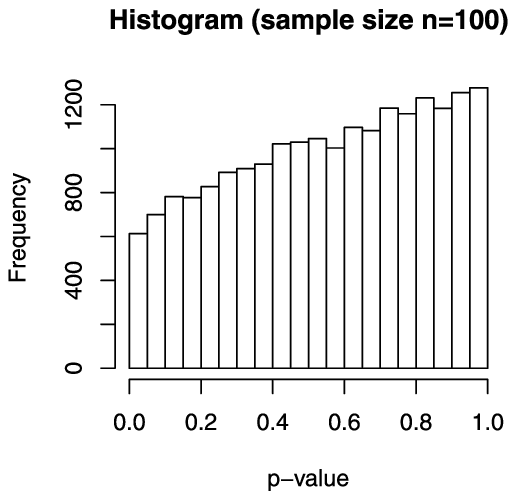}
&\includegraphics{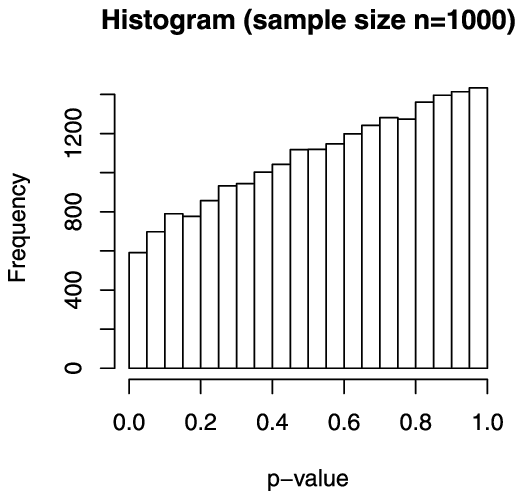}\\
\footnotesize{(a)}&\footnotesize{(b)}\\[6pt]

\includegraphics{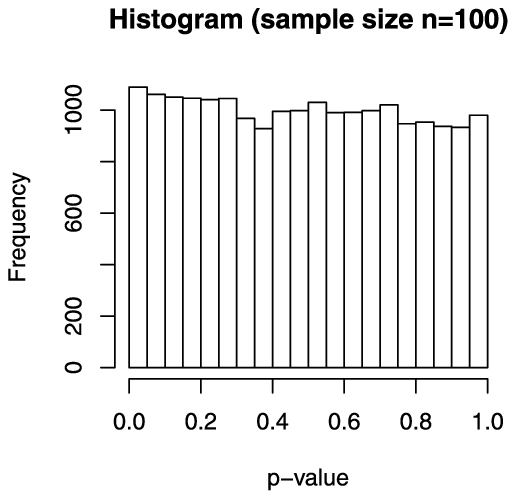}
&\includegraphics{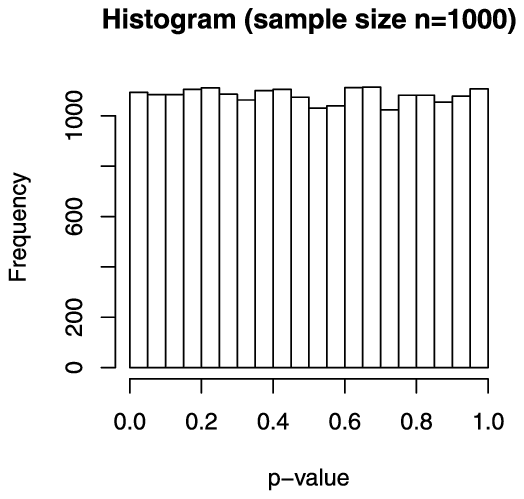}\\
\footnotesize{(c)}&\footnotesize{(d)}
\end{tabular}
\caption{Histograms of $p$-values for a likelihood ratio test.}
\label{fig:histo}
\end{figure}

\begin{example}
\label{ex:simulation}
Consider the graph $G=(V,D,B)$ from Figure~\ref{fig:graph-sim}. Let
$\Lambda=(\lambda_{ij})$ be the matrix in $\mathbb{R}^D$ with
\[
\lambda_{12}=3, \qquad    \lambda_{23} = -\tfrac{1}{2}, \qquad
\lambda_{34}=\lambda_{45}=1.
\]
Let $\Omega=(\omega_{ij})$ be the matrix in $\mathit{PD}(B)$ with all
diagonal entries equal to 2 and
\[
\omega_{14}=\omega_{15}=\omega_{24}=\omega_{35}=1.
\]
It can be shown that at the specified point $(\Lambda,\Omega)$ the map
$\phi_G$ is not injective and the image of $\phi_G$ has a singularity.
%% \[
%% \Lambda=
%% \begin{pmatrix}
%% 0&3&0&0&0\\
%% 0&0&-1/2&0&0\\
%% 0&0&0&1&0\\
%% 0&0&0&0&1\\
%% 0&0&0&0&0
%% \end{pmatrix},
%% \Omega=
%% \begin{pmatrix}
%% 2&0&0&1&1\\
%% 0&2&0&1&0\\
%% 0&0&2&0&1\\
%% 1&1&0&2&0\\
%% 1&0&1&0&2
%% \end{pmatrix}.
%% \]
Suppose we use the likelihood ratio test for testing the model
$\mathcal{M}(G)$ against the saturated alternative given by all
multivariate normal distributions on $\mathbb{R}^m$. The standard
procedure would compare the resulting likelihood ratio statistic to a
chi-square distribution with two degrees of freedom.
Figure~\ref{fig:histo} illustrates the problems with this procedure.
What is plotted are histograms of $p$-values obtained from the chi-square
approximation. Each histogram is based on simulation of 20,000 samples
of size $n=100$ or $n=1000$. The samples underlying the two histograms
in Figure~\ref{fig:histo}(a), (b) are drawn from the multivariate normal
distribution with covariance matrix $\Sigma=\phi_G(\Lambda,\Omega)$ for
the above parameter choices. Many $p$-values being large, it is evident
that the test is too conservative. For comparison, we repeat the
simulations with $\lambda_{23}=1/2$ and all other parameters unchanged.
There is no identifiability failure in this second scenario, the image of
$\phi_G$ is smooth in a neighborhood of the new covariance matrix and, as
shown in Figure~\ref{fig:histo}(c), (d), the expected uniform distribution
for the $p$-values emerges in reasonable approximation. %\qed
\end{example}

Call a directed graph with at least two nodes an arborescence
converging to node $i$ if its edges form a spanning  tree with a
directed path from any node $j\not=i$ to $i$. In other words, $i$ is
the unique sink node. For a mixed graph $G=(V,D,B)$ and a subset of
nodes $A\subset V$, let $D_A=D\cap(A\times A)$ be the set of directed
edges with both endpoints in $A$. Similarly, let $B_A=B\cap(A\times
A)$, and define the mixed subgraph induced by $A$ to be
$G_A=(A,D_A,B_A)$. Our main result provides the following answer to
Question~\ref{que:ident}.

\begin{theorem}
\label{thm:main}
The parametrization $\phi_G$ for an acyclic mixed graph $G=(V,D,B)$
fails to be injective if and only if there is an induced subgraph
$G_A$, $A\subseteq V$, whose directed part $(A,D_A)$ contains a
converging arborescence and whose bidirected part $(A,B_A)$ is
connected. If $\phi_G$ is injective, then its inverse is a rational
map.
\end{theorem}

An acyclic mixed graph $G=(V,D,B)$ is simple if there is at most one
edge between any pair of nodes, that is, if $D\cap B=\varnothing$.
Theorem~\ref{thm:main} states in particular that only simple acyclic
mixed graphs may have an injective parametrization. Indeed, two edges
$i\leftrightarrow j$ and $i\to j$, respectively, connect and yield an arborescence
in the subgraph $G_{\{i,j\}}$.

\begin{corollary}
If the acyclic mixed graph $G$ has at most three nodes, then $\phi_G$
is injective if and only if $G$ is simple. There are exactly two
unlabeled simple acyclic mixed graphs on four nodes with $\phi_G$
not injective.
\end{corollary}

\begin{figure}[b]

\includegraphics{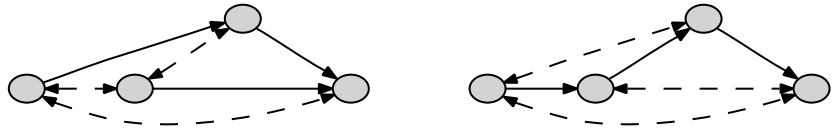}

\caption{The two unlabeled graphs on four nodes with noninjective
parametrization.}
\label{fig:4nodes}
\end{figure}

\begin{pf}
An arborescence involving three nodes contains two edges. The
bidirected part of a simple mixed graph can only be connected if
there are two further edges. However, a simple graph with three
nodes has at most three edges. The two examples on four nodes are
shown in Figure~\ref{fig:4nodes}.
\end{pf}

%
%f3 ###

A possibly cyclic mixed graph $G=(V,D,B)$ is simple if there is at
most one edge between any pair of nodes, that is, if $D\cap
B=\varnothing$ and the presence of an edge $i\to j$ in $D$ implies the
absence of $j\to i$. As shown in the next lemma, it is easy to give a
direct proof of the fact that only simple graphs can have an injective
parametrization. The lemma also clarifies that noninjectivity can be
recognized in subgraphs, which is a fact that is important for later
proofs.

\begin{lemma}
\label{lem:simple}
Suppose the map $\phi_G$ given by a mixed graph $G$ is injective.
Then $G$ is simple, and $\phi_H$ is injective for any
(not necessarily induced\hspace*{1.5pt}) subgraph $H$ of~$G$.
\end{lemma}

\begin{pf}
If $H=(V',D',B')$ is a subgraph of $G=(V,D,B)$, that is,
$V'\subseteq V$, $D'\subseteq D$ and $B'\subseteq B$, then $\phi_H$
is injective if and only if $\phi_G$ is injective at points that
have all parameters $\lambda_{ij}$ and $\omega_{ij}$ zero for edges
$(i,j)\in D\setminus D'$ or $(i,j)\in B\setminus B'$. If $G$ is not
simple, then there exist two distinct indices $i,j$ for which the
graph contains at least two of the three possible edges $i\to j$,
$j\to i$  and $i\leftrightarrow j$. If $V=\{i,j\}$, then $\phi_G$ is not
injective because it maps the at least 4-dimensional set
$\mathbb{R}^D_{\mathrm{reg}}\times\mathit{PD}(B)$ to the 3-dimensional
cone of
positive definite $2\times2$ matrices. If $ |V|>2$, then the claim
follows by passing to the subgraph induced by $\{i,j\}$.
\end{pf}

The remainder of the paper is organized as follows.
Section~\ref{sec:prior-work} reviews the connection of our work to the
existing literature on identifiability of structural equation
models.
%% In particular, we contrast global identifiability with the
%% weaker notion of generic identifiability.
Section~\ref{sec:stepwise-inversion} lays out the natural stepwise
approach to inversion of the parametrization $\phi_G$ in the case
where the underlying graph is acyclic. Necessity and sufficiency of
the graphical condition from our main Theorem~\ref{thm:main} are
proven in Sections~\ref{sec:necessity} and \ref{sec:sufficiency},
respectively. In Section~\ref{sec:lemma_proofs}, we collect three
lemmas used in the proof of sufficiency. Theorem~\ref{thm:acyclic}
about directed cycles is proven in Section~\ref{sec:cyclic-models}.
Concluding remarks are given in Section~\ref{sec:conclusion}.

%s2 ###
\section{Prior work}
\label{sec:prior-work}

Identifiability properties of structural equation models are a topic
with a long history. A review of classical conditions, which do not
take into account the finer graphical structure considered here, can
be found, for instance, in~the monograph \cite{bollen1989}. A more
recent sufficient condition for global identifiability of the linear
structural equation models from Definition~\ref{def:sem} is due to
\cite{mcdonald2002,richardson2002}. It requires the presence of a
bidirected edge $i\leftrightarrow j$ to imply the absence of directed
paths from
$j$ to $i$ (and from $i$ to $j$). Following \cite{richardson2002}, we
call an acyclic mixed graph with this property \textit{ancestral}. It
is clear that an ancestral mixed graph is simple. We revisit the
result about ancestral graphs in Corollary~\ref{cor:ancestral} below.

Other recent work, such as \cite{brito2002}, considers a weaker
identifiability requirement for the model $\mathcal{M}(G)$ associated with
a mixed graph $G=(V,D,B)$. For a pair of matrices
$\Lambda_0\in\mathbb{R}^D_{\mathrm{reg}}$ and $\Omega_0\in\mathit{PD}(B)$,
define the fiber
%
%e2.1 ###
\begin{equation}
\label{eq:fiber} \qquad
\mathcal{F}(\Lambda_0,\Omega_0) =  \{ (\Lambda,\Omega) \dvtx
\phi_G(\Lambda,\Omega)=\phi_G(\Lambda_0,\Omega_0) ,
\Lambda\in\mathbb{R}^D_{\mathrm{reg}},  \Omega\in\mathit{PD}(B)
 \}.
\end{equation}
The map $\phi_G$ is injective if and only if all its fibers contain
only a
single point. If it holds instead that for generic choices of
$\Lambda\in\mathbb{R}^D_{\mathrm{reg}}$ and $\Omega\in\mathit
{PD}(B)$, the fiber
$\mathcal{F}(\Lambda,\Omega)$ contains only the single point
$(\Lambda,\Omega)$, then we say that the map $\phi_G$ is \textit{generically
injective} and the model $\mathcal{M}(G)$ is \textit{generically
identifiable}. Requiring a condition to hold for generic points means
that the points at which the condition fails form a lower-dimensional
algebraic subset. In particular, the condition holds for almost every
point (in Lebesgue measure), and some authors thus also speak of an almost
everywhere identifiable model; compare the lemma in \cite{okamoto1973}.
When the substantive interest is in all parameters of a model, generic
identifiability constitutes a minimal requirement. However, generically
but not globally identifiable models can have nonsmooth parameter spaces
and thus present difficulties for statistical inference; recall
Example~\ref{ex:simulation} that treats a generically identifiable model.

The main theorem of \cite{brito2002}, which we reprove in
Corollary~\ref{cor:brito}, states that $\phi_G$ is generically injective
for every simple acyclic mixed graph $G$. The graph being simple and
acyclic, however, is far from necessary for generic injectivity of
$\phi_G$. A classical counterexample is the instrumental variable model
based on the graph with edges $1\to2\to3$ and $2\leftrightarrow3$.
Cyclic models
may also be generically identifiable; for instance, see Example~3.6 in
\cite{drtonlrt}. For recent work on the topic, see \cite{tian2009} and
references therein. To our knowledge, characterizing the mixed graphs $G$
with generically injective parametrization $\phi_G$ remains an open
problem.

The linear structural equation models $\mathcal{M}(G)$ considered in this
paper are closely related to latent variable models known as semi-Markovian
causal models. These nonparametric models are obtained by subdividing the
bidirected edges, that is, each edge $i\leftrightarrow j$ is replaced
by two directed
edges $i\leftarrow u_{ij} \to j$, where $u_{ij}$ is a new node. Each node
$u_{ij}$ added to the vertex set corresponds to a latent variable; compare
also \cite{richardson2002,pearl2009,wermuth2010}. Using results from
\cite{tian2002}, the work of \cite{shpitser2006} gives graphical
conditions for when (univariate or multivariate) intervention distributions
in acyclic semi-Markovian causal models are identified. This work is based
on manipulating recursive density factorizations involving latent
variables. If $G$ is an acyclic mixed graph and the structural equation
model $\mathcal{M}(G)$ is contained in the semi-Markovian model for $G$,
then $\mathcal{M}(G)$ is globally identified provided that in the
semi-Markovian model we can identify, for every node $i$, the univariate
intervention distribution for $i$ and intervention set $\pa(i)$; see also
Chapter~6 in \cite{tian2002}.

For an acyclic mixed graph $G=(V,D,B)$, we may define a Gaussian model
$\mathcal{M}'(G)$ by assuming that both the observed and the latent
variables in the semi-Markovian model for $G$ have a joint multivariate
normal distribution. This creates an explicit connection to linear
structural equation models, and it is indeed possible that
$\mathcal{M}'(G)=\mathcal{M}(G)$. For instance, if there are no directed
edges $(D=\varnothing)$, then $\mathcal{M}'(G)=\mathcal{M}(G)$ if and
only if
the bidirected part $(V,B)$ is a forest of trees; see Corollary~3.4 in
\cite{drtonyu2010}. If $D=\varnothing$ and $(V,B)$ is not a forest of
trees, then $\mathcal{M}(G)$ is strictly larger than $\mathcal{M}'(G)$.
Therefore, other nonnormal constructions would be required in order for
the theorems in \cite{shpitser2006} to furnish sufficient conditions for
global identifiability of linear structural equation models. We are
unaware, however, of literature providing a connection between
semi-Markovian causal models and the linear structural equation models from
Definition~\ref{def:sem} when non-Gaussian distributions are assumed for
the latent variables.

Finally, the existing counterexamples to identifiability of semi-Markovian
models involve binary variables and thus cannot be used to prove necessity
of an identifiability condition for the Gaussian models $\mathcal{M}(G)$.
However, despite this fact and the difficulties in relating the models
$\mathcal{M}(G)$ to semi-Markovian models, our graphical condition from
Theorem~\ref{thm:main}, which we first found by experimentation with
computer algebra software, coincides with that of \cite{shpitser2006}; the
term ``$y$-rooted C-tree'' is used there to refer to a mixed graph whose
directed part is an arborescence converging to node $y$ and whose
bidirected part is a tree. A reader familiar with the work in
\cite{tian2002} will also recognize similarities between the higher-level
structure of the proofs given there and those in
Section~\ref{sec:sufficiency} of this paper.

%s3 ###
\section{Stepwise inversion}
\label{sec:stepwise-inversion}

Throughout this section, suppose that $G\!=\!(\!V,D,B\!)$ is an acyclic mixed graph
with vertex set $V=[m]$. The map $\phi_G$ is injective if all its fibers
contain only a single point; recall the definition of a fiber in
(\ref{eq:fiber}). Let $\Sigma=\phi_G(\Lambda_0,\Omega_0)$ for two matrices
$\Lambda_0\in\mathbb{R}^D$ and $\Omega_0\in\mathit{PD}(B)$. This section
describes how to find points $(\Lambda,\Omega)$ in the fiber
$\mathcal{F}(\Lambda_0,\Omega_0)$. In particular, we show in
Lemma~\ref{lem:overall-ident} that an algebraic criterion can be used to
decide whether the map $\phi_G$ is injective. The lemma is proven
after we
describe a~natural inversion approach that uses the acyclic structure of
the graph $G$ in a stepwise manner. We remark that this stepwise inversion
is closely related to the idea of pseudo-variable regression used in the
iterative conditional fitting algorithm of \cite{drtoneichler2009}.
%% at a given point $(\Lambda,\Omega)$.

For each $i \leq m-1$, let $P(i) = \operatorname{pa}(i+1)$ be the parents of
node $i+1$, and $S(i) = \{ j \leq i \dvtx  j \leftrightarrow i+1 \in B \}$
the siblings
of $i+1$. (In other related work, the nodes incident to a bidirected
edge $i\leftrightarrow j$ have also been called ``spouses'' of each
other but we
find ``siblings'' to be natural terminology given that a~common parent
to the two nodes is introduced when subdividing the edge as discussed
in Section~\ref{sec:prior-work}.)

\begin{lemma}
\label{lem:overall-ident}
Suppose $G=(V,D,B)$ is an acyclic mixed graph with its nodes labeled
in a topological order. Then the parametrization $\phi_G$ is injective
%% fiber $\mathcal{F}(\Lambda,\Omega)$
%% contains a single point
if and only if the rank condition
\[
\rank \bigl(  \Omega_{[i]\setminus
S(i),[i]}(I-\Lambda)^{-1}_{[i], P(i)}
 \bigr) = | P(i)|
\]
holds for all nodes $i=1,\ldots ,m-1$ and all pairs
$\Lambda\in\mathbb{R}^D$ and $\Omega\in\mathit{PD}(B)$.
\end{lemma}

\begin{remark}
In this paper, matrix inversion is always given higher priority than an
operation of forming a submatrix. For any invertible matrix $M$ and
index sets $A,B$, the matrix $M^{-1}_{A,B}=(M^{-1})_{A,B}$ is thus the
$A\times B$ submatrix of the inverse of $M$.
\end{remark}

Computing points $(\Lambda,\Omega)$ in the fiber
$\mathcal{F}(\Lambda_0,\Omega_0)$ means solving the polynomial equation
system given by the matrix equation
%
%e3.1 ###
\begin{equation}
\label{eq:sigma}
\Sigma=(I-\Lambda)^{-T}\Omega(I-\Lambda)^{-1}.
\end{equation}
For topologically ordered nodes, (\ref{eq:sigma}) implies that
$\sigma_{11}=\omega_{11}$ and that the first column in the strictly
upper-triangular matrix $\Lambda$ contains only zeros. Hence, these are
uniquely determined for all matrices in the fiber.

Let $i\ge1$, and assume that we know the $[i]\times[i]$ submatrices of
$\Lambda$ and $\Omega$ of a solution to   equation   (\ref{eq:sigma}).
Partition off the $(i+1)$st row and column of the submatrices
\[
(I-\Lambda)_{[i+1], [i+1]} =
\pmatrix{
\Gamma& -\lambda\cr
0 & 1
}
,
\qquad
\Omega_{[i+1], [i+1]} =
\pmatrix{
\Psi& \omega\cr
\omega^T & \omega_{i+1,i+1}
}
.
\]
The matrices $\Gamma$ and $\Psi$ are known, $\lambda_{[i]\setminus P(i)}=0$
and $\omega_{[i]\setminus S(i)}=0$. The inverse of
$I-\Lambda$ can be written as a block matrix as
%
%e3.2 ###
\begin{equation}
\label{eq:inverse-L}
(I-\Lambda)^{-1}_{[i+1], [i+1]}
%% := [(I-\Lambda)^{-1}]_{[i+1], [i+1]}
=
\pmatrix{
\Gamma^{-1} & \Gamma^{-1}\lambda\cr
0 & 1
}
.
\end{equation}

In this notation, the part of equation (\ref{eq:sigma}) that
pertains to the $[i+1]\times[i+1]$ submatrix of $\Sigma$ is
\begin{eqnarray*}
&&\pmatrix{
\Sigma_{[i], [i]} & \Sigma_{[i],\{i+1\}}\cr
&\sigma_{i+1,i+1}
}\\
&& \qquad =
\pmatrix{
\Gamma^{-T}\Psi\Gamma^{-1} & \Gamma^{-T}\Psi\Gamma^{-1} \lambda+
\Gamma^{-T} \omega\cr
& \omega_{i+1,i+1} + \lambda^T \Gamma^{-T}\Psi\Gamma^{-1}\lambda+
2\omega^T\Gamma^{-1}\lambda
}
,
\end{eqnarray*}
where only the upper-triangular parts of the symmetric matrices are shown.
Hence, given the values of $\Gamma$ and $\Psi$, the choice of
$\lambda$
and $\omega$ is unique if and only if the equation
%
%e3.3 ###
\begin{equation}
\label{eq:identify1}
\Sigma_{[i],\{i+1\}} = \Gamma^{-T}\Psi\Gamma^{-1}\cdot\lambda+
\Gamma^{-T}\cdot\omega
\end{equation}
has a unique solution. Clearly, any feasible choice of a solution
$(\lambda,\omega)$ to the equation in
(\ref{eq:identify1}) leads to a unique solution $\omega_{i+1,i+1}$
via the
equation
%
%e3.4 ###
\begin{equation}
\label{eq:identify2}
\sigma_{i+1,i+1} = \omega_{i+1,i+1} + \lambda^T
\Gamma^{-T}\Psi\Gamma^{-1}\lambda+
2\omega^T\Gamma^{-1}\lambda.
\end{equation}
Since $\lambda_{[i]\setminus P(i)}=0$ and $\omega_{[i]\setminus S(i)}=0$,
equation (\ref{eq:identify1}) can be rewritten as
\[
\Sigma_{[i],\{i+1\}} = \bigl(\Gamma^{-T}\Psi\Gamma^{-1}_{[i],
P(i)}\bigr)\cdot
\lambda_{ P(i)} + \bigl(\Gamma^{-1}_{ S(i),[i]}\bigr)^{T}\cdot\omega_{ S(i)}.
\]
It has a unique solution if and only if the matrix
\[
\left[\vphantom{\bigl(\Gamma^{-1}_{ S(i),[i]}\bigr)^{T}}
\matrix{
\Gamma^{-T}\Psi\Gamma^{-1}_{[i], P(i)} & \bigl(\Gamma^{-1}_{ S(i),[i]}\bigr)^{T}
} \right]
\]
has full column rank $| P(i)|+| S(i)|$. The matrix $\Gamma$ is
invertible because it~is upper-triangular with ones along the diagonal.
Thus, the condition is equivalent to
\[
\Gamma^T
 \left[\vphantom{\bigl(\Gamma^{-1}_{ S(i),[i]}\bigr)^{T}}
  \matrix{
\Gamma^{-T}\Psi\Gamma^{-1}_{[i], P(i)} & \bigl(\Gamma^{-1}_{ S(i),[i]}\bigr)^{T}
} \right]
=
\left[\vphantom{\Psi\Gamma^{-1}_{[i], P(i)}} \matrix{
\Psi\Gamma^{-1}_{[i], P(i)} & I_{[i], S(i)}
} \right]
\]
having full column rank. The second block is part of an identity
matrix. We deduce that the condition is equivalent to requiring
that $\Psi_{[i]\setminus S(i),[i]}\Gamma^{-1}_{[i], P(i)}$, the\vspace*{1pt}
submatrix obtained by removing the rows and columns with
index in $ S(i)$, has rank $| P(i)|$. Note that
\[
\Psi_{[i]\setminus S(i),[i]}\Gamma^{-1}_{[i], P(i)}=\Omega
_{[i]\setminus
S(i),[i]}(I-\Lambda)^{-1}_{[i], P(i)}
\]
is the matrix appearing in Lemma~\ref{lem:overall-ident}.

\begin{pf*}{Proof of Lemma~\ref{lem:overall-ident}}
Consider a feasible pair $(\Lambda,\Omega)$. If the rank condition for
this pair holds for all nodes $i=1,\ldots ,m-1$, then it follows from the
stepwise inversion procedure described above that the fiber
$\mathcal{F}(\Lambda,\Omega)$ contains only the single point
$(\Lambda,\Omega)$. Therefore, the rank condition holding for all nodes
and all matrix pairs implies that all fibers are singletons, or in other
words, that the map $\phi_G$ is injective.

Conversely, assume that the rank condition fails for some node \mbox{$i\le m-1$}
and matrix pair $(\Lambda,\Omega)$. If $i=m-1$, then the considered
fiber $\mathcal{F}(\Lambda,\Omega)$ is positive-dimensional, and
$\phi_G$
not injective. If $i<m-1$, then it follows analogously that the
parametrization $\phi_H$ for the induced subgraph $H=G_{[i+1]}$ is not
injective. By Lemma~\ref{lem:simple}, $\phi_G$ cannot be injective
either.
\end{pf*}

If the rank condition in Lemma~\ref{lem:overall-ident} holds at a
particular pair $(\Lambda,\Omega)$, then the fiber
$\mathcal{F}(\Lambda,\Omega)$ contains only the pair $(\Lambda
,\Omega)$.
However, the converse is false in general, that is, failure of the rank
condition at a particular pair $(\Lambda,\Omega)$ and vertex $i < m$
need not imply that the
fiber $\mathcal{F}(\Lambda,\Omega)$ contains more than one point.
This may
occur even for a simple acyclic mixed graph.

\begin{figure}

\includegraphics{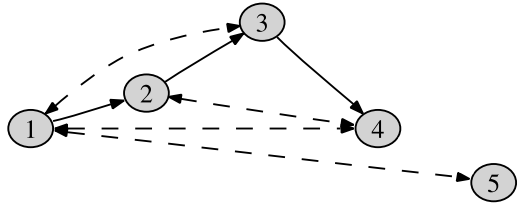}

\caption{Graph with noninjective parametrization (see
Example~\protect\ref{ex:not-extend}).}
\label{fig:6nodes}
\end{figure}

\begin{example}
\label{ex:not-extend}
Consider the graph in Figure~\ref{fig:6nodes}, set
$\lambda_{12}=\lambda_{23}=\lambda_{34}=1$, and choose the positive
definite matrix
\[
\Omega=
\pmatrix{
2&0&-1&-1&-1\cr
0&1&0&-1&0\cr
-1&0&1&0&0\cr
-1&-1&0&3&0\cr
-1&0&0&0&3
}
.
\]
The rank condition for this pair $(\Lambda,\Omega)$ fails at node $i=3$.
Nevertheless, the fiber $\mathcal{F}(\Lambda,\Omega)$ is equal to
$\{(\Lambda,\Omega)\}$. If we set $\omega_{15}=0$, however, then
$\mathcal{F}(\Lambda,\Omega)$ becomes one-dimensional. Using terminology
from econometrics/causality, the variable corresponding to node 5 behaves
like an ``instrument;'' compare, for instance, \cite{pearl2009}.
\end{example}

%
%f4 ###

Lemma \ref{lem:overall-ident} allows us to give simple proofs of two
established results in the graphical models literature. The proof of
Corollary~\ref{cor:ancestral} emphasizes the special structure
exhibited by
ancestral graphs. The proof of Corollary~\ref{cor:brito} demonstrates that
the identity matrix always has a singleton as a fiber under the
parametrization associated with a simple acyclic mixed graph.

\begin{corollary}
\label{cor:ancestral}
If the acyclic mixed graph $G$ is ancestral then the parame\-trization
$\phi_G$ is injective.
\end{corollary}

\begin{pf}
Recall that if $G=(V,D,B)$ is ancestral and $i\leftrightarrow j$ is a
bidirected edge
in $G$, then there is no directed path from $i$ to $j$ or $j$ to $i$.
Suppose $V=[m]$ is topologically ordered, and let $i$ be some node
smaller than $m$. Pick a node $j\in S(i)$. Then there may not exist a
directed path from $j$ to a~node in $P(i)$. It follows that
\[
\Omega_{[i]\setminus
S(i),[i]}(I-\Lambda)^{-1}_{[i], P(i)} =
\Omega_{[i]\setminus
S(i),[i]\setminus
S(i)}(I-\Lambda)^{-1}_{[i]\setminus
S(i), P(i)}.
\]
The latter matrix is the product of a principal and thus positive
definite submatrix of $\Omega$ and a matrix that contains the
$P(i)\times
P(i)$ identity matrix. It follows that this product has full column rank
$|P(i)|$ for all feasible pairs $(\Lambda,\Omega)$ and all nodes
$i\le
m-1$. By Lemma~\ref{lem:overall-ident}, $\phi_G$ is injective.
\end{pf}

If the acyclic mixed graph $G$ is simple, then $P(i)\subseteq
[i]\setminus
S(i)$ for all nodes $i\le m-1$. Hence, the matrix product appearing in the
rank condition always has at least as many rows as columns. The next
generic identifiability result follows immediately; recall the definitions
in Section~\ref{sec:prior-work}.

\begin{corollary}
\label{cor:brito}
If $G=(V,D,B)$ is a simple acyclic mixed graph, then~the map $\phi_G$ is
generically injective.
\end{corollary}

\begin{pf}
We need to show that for generic choices of $\Lambda\in\mathbb{R}^D$ and
$\Omega\in\mathit{PD}(B)$, the fiber $\mathcal{F}(\Lambda,\Omega
)$ is
equal to the singleton $\{(\Lambda,\Omega)\}$. Set $\Lambda=0$ and
choose $\Omega$ to be the identity matrix. Then each of the matrix
products
%
%e3.5 ###
\begin{equation}
\label{eq:brito}
\Omega_{[i]\setminus
S(i),[i]}(I-\Lambda)^{-1}_{[i], P(i)}, \qquad i=1,\ldots , m-1,
\end{equation}
has the identity matrix as $P(i)\times P(i)$ submatrix. The rank
condition from Lemma~\ref{lem:overall-ident} thus holds for all $i\le
m-1$. Since the matrices in (\ref{eq:brito}) have polynomial entries,
existence of a single pair $(\Lambda,\Omega)$ at which the $m-1$ matrices
in (\ref{eq:brito}) have full column rank implies that the set of pairs
$(\Lambda,\Omega)$ for which at least one of the matrices fails to have
full column rank is a~lower-dimensional algebraic set; compare
\cite{cox2007}, Chapter~9, for background on such algebraic arguments.
\end{pf}

In order to prepare for arguments turning the algebraic condition from
Lemma~\ref{lem:overall-ident} into a graphical one, we detail the structure
of the inverse $(I-\Lambda)^{-1}$ for a matrix
$\Lambda=(\lambda_{ij})\in\mathbb{R}^D$. Let $\mathcal{P}(i,j)$
denote the
set of directed paths from $i$ to $j$ in the considered acyclic graph.

\begin{lemma}
\label{lem:inverse-paths}
The entries of the inverse $(I-\Lambda)^{-1}$ are
\[
\label{eq:inverse-paths}
(I-\Lambda)^{-1}_{ij} = \sum_{\pi\in\mathcal{P}(i,j)} \prod_{k\to
l\in\pi}
\lambda_{kl}, \qquad i,j\in[m].
\]
\end{lemma}

\begin{pf}
This well-known fact can be shown by induction on the matrix size $m$ and
using the partitioning in (\ref{eq:inverse-L}) under a topological
ordering of the nodes.
\end{pf}

Note that adopting the usual definition that takes an empty sum to be zero
and an empty product to be one, the formula in
Lemma~\ref{lem:inverse-paths} states that $(I-\Lambda)^{-1}_{ij}=0$ if
$i\not=j$ and $\mathcal{P}(i,j)=\varnothing$, and it states that
$(I-\Lambda)^{-1}_{ii}=1$ because $\mathcal{P}(i,i)$ contains only a
trivial path without edges.

%s4 ###
\section{Necessity of the graphical condition for identifiability}
\label{sec:necessity}

We now prove that the graphical condition in Theorem~\ref{thm:main}, which
states that there be no induced subgraph whose directed part contains a
converging arborescence and whose bidirected part is connected, is
necessary for the parametrization $\phi_G$ to be injective. By
Lemma~\ref{lem:simple}, it suffices to consider an acyclic mixed graph
whose directed part is a converging arborescence and whose bidirected part
is a spanning tree. In light of Lemma~\ref{lem:overall-ident}, the
necessity of the graphical condition in Theorem~\ref{thm:main} then follows
from the following result.

\begin{proposition}
\label{prop:necessity}
Let $G=(V,D,B)$ be an acyclic mixed graph with topologically ordered
vertex set $V=[m+1]$. If $(V,D)$ is an arborescence converging to $m+1$
and $(V,B)$ is a spanning tree, then there exists a pair of matrices
$\Lambda\in\mathbb{R}^D$ and $\Omega\in\mathit{PD}(B)$ with
\[
\kernel \bigl(  \Omega_{[m]\setminus S(m),[m]}(I-\Lambda)^{-1}_{[m],
P(m)}    \bigr) \not=\{0\}.
\]
\end{proposition}

Let $\mathcal{L}(\Lambda)\subseteq\mathbb{R}^{m}$ be the column
span of
$(I-\Lambda)^{-1}_{[m], P(m)}$. We formulate a~first lemma that we will
use to prove Proposition~\ref{prop:necessity}.

\begin{lemma}
\label{lem:union-spans}
If $V=[m+1]$ and $(V,D)$ is an arborescence converging to node $m+1$,
then the union of the linear spaces $\mathcal{L}(\Lambda)$ for all
$\Lambda\in\mathbb{R}^D$ contains the set
$(\mathbb{R}^*)^m=(\mathbb{R}\setminus\{0\})^m$ of vectors with all
coordinates nonzero.
\end{lemma}

\begin{pf}
In the arborescence, there is a unique path $\pi(i)$ from any vertex
$i\in[m]\setminus P(m)$ to the sink node $m+1$. Let $k(i)$ be the unique
node in $P(m)$ that lies on this path. Let $\Lambda\in\mathbb{R}^D$
and $\alpha\in\mathbb{R}^{|P(m)|}$, and define the vector
\[
\beta(\Lambda,\alpha)=(I-\Lambda)^{-1}_{[m], P(m)}\alpha\in
\mathbb{R}^m.
\]
Since the principal submatrix $(I-\Lambda)^{-1}_{P(m), P(m)}$ is an
identity matrix (because the directed graph is a converging
arborescence), $\beta(\Lambda,\alpha)_i=\alpha_i$ for all $i\in
P(m)$. For
$i\in[m]\setminus P(m)$, we use Lemma~\ref{lem:inverse-paths} to obtain
%
%e4.1 ###
\begin{equation}
\label{eq:torus-recursion}
\beta(\Lambda,\alpha)_i = \alpha_{k(i)}\prod_{j\to l\in\pi(i)}
\lambda_{jl}
= \lambda_{ij} \beta(\Lambda,\alpha)_j,
\end{equation}
where $i\to j\in G$ is the unique edge originating from $i$.

Let $x$ be any vector in $(\mathbb{R}^*)^m$. Our claim states that there
exist a matrix $\Lambda\in\mathbb{R}^D$ and vector $\alpha$ such that
$x=\beta(\Lambda,\alpha)$. Clearly, $\alpha$ has to be equal to the
subvector $x_{P(m)}$. The associated unique choice of $\Lambda$ is
obtained by recursively solving for the entries $\lambda_{ij}$ using the
relationship in (\ref{eq:torus-recursion}).
\end{pf}

Let $R(m)=[m]\setminus S(m)$ be the ``rest'' of the nodes. We are
left with the problem of finding a matrix $\Omega\in\mathit{PD}(B)$
for which some vector in $(\mathbb{R}^*)^m$ lies in the kernel of the
submatrix
\[
\Omega_{R(m),[m]} =
\left[\vphantom{\Omega_{R(m),R(m)}} \matrix{
\Omega_{R(m),R(m)} &
\Omega_{R(m),S(m)}
} \right]
.
\]
%
%We note that if $\Omega_{R(m),R(m)}$ and $\Omega_{R(m),S(m)}$ have the
%required zeros and $\Omega_{R(m),R(m)}$ is positive definite,
%then there is a completion to a positive definite matrix
%$\Omega\in\mathit{PD}(B)$.
Proposition~\ref{prop:necessity} now follows by
combining Lemma~\ref{lem:union-spans} with the next result.

\begin{lemma}
\label{lem:laplacian}
If $(V,B)$ is a tree on $V=[m+1]$, then there exists a matrix
$\Omega\in\mathit{PD}(B)$ such that the vector $\mathbf{1}=(1,\ldots ,1)^T$
is in the kernel of the submatrix $\Omega_{R(m),[m]}$.
\end{lemma}

\begin{pf}
Let $T$ be the set of all nodes in $R(m)$ that are connected to some node
in $S(m)$ by an edge in $B$. If $\Omega\in\mathit{PD}(B)$, then the
submatrix $\Omega_{R(m),S(m)}$ has only zero entries in rows indexed by
nodes $i\in R(m)\setminus T$. If $i\in T$, then the $i$th row of
$\Omega_{R(m),S(m)}$ has at least one entry that is not constrained to
zero and may take any real value. Hence, we can choose a~matrix
$\Omega_{R(m),S(m)}$ that has row sum
%
%e4.2 ###
\begin{equation}
\label{eq:Omega-part1}
\sum_{j\in S(m)}\omega_{ij} =
\cases{\displaystyle
-1 ,&\quad  if   $i\in T$,\cr\displaystyle
0 ,&\quad  if   $i\in R(m)\setminus T$.
}
\end{equation}

Let $H = (R(m),B_{R(m)})$ be the induced subgraph of $G$ on vertex set $R(m)$.
The Laplacian of $H$, $L(H)=(l_{ij})$,
is the symmetric $R(m)\times R(m)$ matrix whose diagonal entries are the
degrees of the nodes in $H$ and whose off-diagonal entries $l_{ij}$ are equal
to $-$1 if $i\leftrightarrow j$ is an edge in $H$ and 0 otherwise. The
Laplacian is
well known to be positive semidefinite with all row sums zero. For a
subset $C\subset[m]$, let $\mathbf{1}_C\in\mathbb{R}^m$ be the vector
with entries equal to one at indices in $C$ and zero elsewhere. The
kernel of $L(H)$ is the direct sum of the linear spaces spanned by the
vectors $\mathbf{1}_C$ for the connected components $C$ of the graph $H$;
compare \cite{chung1997}, Chapter~1.

Let $D_T=(d_{ij})$ be the diagonal matrix that has diagonal entry
$d_{ii}=1$ if $i\in T$ and $d_{ii}=0$ otherwise. Both $L(H)$ and
$D_T$ are positive semidefinite matrices and thus the kernel of
$L(H) + D_T$ is equal to $\ker L(H) \cap\ker D_T$. Since $(V,B)$
is a connected graph, each connected component of $H$ contains a
node in $T$. Therefore, none of the vectors $\mathbf{1}_C$ are in
the kernel of $D_T$, where $C$ ranges over all connected components
of $H$. This implies that the $\ker(L(H) + D_T) = \{0\}$, and
hence this matrix is positive definite.

Let $\Omega$ be any matrix in
$\mathit{PD}(B)$ whose submatrix $\Omega_{R(m),S(m)}$ satisfies
(\ref{eq:Omega-part1}) and whose principal submatrix $\Omega_{R(m),R(m)}$
is the positive definite matrix $L(H)+D_T$. The matrix
$\Omega\in\mathit{PD}(B)$ has the desired property because
\[
\Omega_{R(m),[m]}\mathbf{1}=
\bigl(L(H)+D_T\bigr)\mathbf{1} + \Omega_{R(m),S(m)}\mathbf{1} =
\mathbf{1}_{T} - \mathbf{1}_T = 0.
\]
Such matrices exist because we can choose $\Omega_{S(m), S(m)}$ to
be, for instance, a diagonal matrix with very large diagonal
entries. Principal minors of $\Omega$ that are not submatrices of
$\Omega_{R(m), R(m)}$ will be dominated by these diagonal entries
and hence be positive. All other principal minors are positive
since $\Omega_{R(m),R(m)} = L(H) + D_T$ was shown to be positive
definite.
\end{pf}

%s5 ###
\section{Sufficiency of the graphical condition for identifiability}
\label{sec:sufficiency}

In this section, we prove that the graphical condition in
Theorem~\ref{thm:main}, which requires an acyclic mixed graph $G$ to
have no
induced subgraph whose directed part contains a converging arborescence and
whose bidirected part is connected, is sufficient for the parametrization
$\phi_G$ to be injective. Proposition~\ref{prop:sufficiency} below shows
that if $\phi_G$ is not injective and $G$ does not contain an induced
subgraph with both a converging arborescence and a bidirected spanning
tree, then there is a subgraph $G'$ with fewer nodes such that $\phi_{G'}$
still fails to be injective. The sufficiency of the graphical condition
then follows immediately. To see this, note that a graph $G$ with
noninjective parametrization $\phi_G$ must contain some minimal induced
subgraph $G'$ with noninjective $\phi_{G'}$. Applying the contrapositive
of Proposition~\ref{prop:sufficiency} to $G'$, we conclude that the
directed part of $G'$ contains a converging arborescence and the bidirected
part of $G'$ is connected.

In preparing for the proof of Proposition~\ref{prop:sufficiency}, we
first treat the case when there is no arborescence; this gives
Proposition~\ref{prop:directed_sufficiency}. The case when there is
no bidirected spanning tree is treated in
Proposition~\ref{prop:bidirected_sufficiency}. In either case, we
reduce a given graph $G=(V,D,B)$ to the subgraph $G_W$ induced by a
subset $W\subsetneq V$. We use the notation $\tilde{\Lambda}$,
$\tilde{\Omega}$, $\tilde{P}(i)$, $\tilde{S}(i)$,
$\tilde{\mathcal{P}}(i,j)$ to denote the counterparts to $\Lambda$,
$\Omega$, $P(i)$, $S(i)$  and $\mathcal{P}(i,j)$, when performing this
reduction of $G$ to $G_W$.\vspace{-2pt}

\begin{proposition}
\label{prop:directed_sufficiency}
Let $G=(V,D,B)$ be an acyclic mixed graph with topologically ordered
vertex set $V=[m+1]$, with some $\Lambda\in\mathbb{R}^D$,
$\Omega\in\mathit{PD}(B)$  and nonzero
$\alpha\in\mathbb{R}^{|P(m)|}$, such that\vspace{-4pt}
\[
\Omega_{[m]\setminus
S(m),[m]}(I-\Lambda)^{-1}_{[m], P(m)}\alpha=0.\vspace{-3pt}
\]
Suppose the directed part of $G$ does not contain an arborescence
converging to $m+1$. Let $A$ be the set of nodes $i\leq m$ with some
path of directed edges from $i$ to $m+1$, and $W=A\cup\{m+1\}$.
Then $W\subsetneq V$ and $\phi_{G_W}$ is not injective.\looseness=-1
\vspace{-2pt}
\end{proposition}

\begin{pf}
Since $G$ does not have a converging arborescence, $A\subsetneq[m]$
and \mbox{$W\subsetneq V$}.

Denote the induced subgraph as $G_W=(W,\tilde D,\tilde B)$. Let
$\tilde{\Lambda}=\Lambda_{W,W}\in\mathbb{R}^{\tilde D}$ and
$\tilde{\Omega}=\Omega_{W,W}\in\mathit{PD}(\tilde B)$. Note that
$P(m)\subseteq A$ by definition, and so $\tilde{P}(m)=P(m)$.
Suppose $j\in P(m)$. Then for each $i\in[m]\setminus A$,
$\mathcal{P}(i,j)=\varnothing$ by definition, and so
$(I-\Lambda)^{-1}_{ij}=0$ by Lemma~\ref{lem:inverse-paths}. For
each $i\in A$, and for any path $i\to v_1\to\cdots \to v_k\to j$ in
$G$, each intermediate vertex $v_1,\ldots ,v_k$ is in $A$ by
definition of $A$ (since there is an edge $j\to m+1$). Therefore,
$\tilde{\mathcal{P}}(i,j)=\mathcal{P}(i,j)$, and it follows that
$(I-\tilde\Lambda)^{-1}_{ij}=(I-\Lambda)^{-1}_{ij}$. In other
words, when the nodes outside of $W$ are removed from $G$, the
remaining entries of $(I-\Lambda)^{-1}$ are unchanged, while the
removed entries in the columns indexed by $P(m)=\tilde P(m)$ are all
zero. We obtain that\vspace{-2pt}
\begin{eqnarray*}
\sum_{i\in A}\tilde{\Omega}_{A\setminus\tilde
S(m),i}(I-\tilde{\Lambda})^{-1}_{i,\tilde P(m)}\alpha
&=& \sum_{i\in A}\Omega_{A\setminus
S(m),i}(I-\Lambda)^{-1}_{i,P(m)}\alpha\\
&=&\sum_{i\in[m]}
\Omega_{A\setminus S(m),i}(I-\Lambda)^{-1}_{i,P(m)}\alpha\\
&=&
\Omega_{A\setminus S(m),[m]}(I-\Lambda)^{-1}_{[m],P(m)}\alpha.
\end{eqnarray*}
By assumption, the last quantity is zero. By
Lemma~\ref{lem:overall-ident}, $\phi_{G_W}$ is not injective.
\end{pf}

We next prove a similar proposition for graphs whose bidirected part is not
connected. The proof uses Lemmas~\ref{lem:sublemma1} and
\ref{lem:sublemma3}, which are derived in Section~\ref{sec:lemma_proofs}.

\begin{proposition}
\label{prop:bidirected_sufficiency}
Let $G=(V,D,B)$ be an acyclic mixed graph with topologically ordered
vertex set $V=[m+1]$, with some $\Lambda\in\mathbb{R}^D$,
$\Omega\in\mathit{PD}(B)$, and nonzero $\alpha\in\mathbb{R}^{|P(m)|}$,
such that
\[
\Omega_{[m]\setminus
S(m),[m]}(I-\Lambda)^{-1}_{[m], P(m)}\alpha=0.
\]
Suppose the bidirected part of $G$ is not connected. Let $A$ be the set
of nodes $i\leq m$ with some path of bidirected edges from $i$ to $m+1$,
and $W=A\cup\{m+1\}$. Then $W\subsetneq V$ and $\phi_{G_W}$ is not
injective.
\end{proposition}

\begin{pf}
Since the bidirected part is not connected, $A\subsetneq[m]$ and
$W\subsetneq V$.

Denote the induced subgraph as $G_W=(W,\tilde D,\tilde B)$. Let
$\tilde{\Lambda}=\Lambda_{W,W}\in\mathbb{R}^{\tilde D}$ and
$\tilde{\Omega}=\Omega_{W,W}\in\mathit{PD}(\tilde B)$. If $i\in S(m)$,
then it holds trivially that $i\in A$ and thus $\tilde{S}(m)=S(m)$. By
Lemma~\ref{lem:sublemma3} below,
\begin{eqnarray*}
\tilde{\Omega}_{A\setminus
\tilde
S(m),A}(I-\tilde\Lambda)^{-1}_{A,\tilde{P}(m)}\alpha_{\tilde{P}(m)}
&=&\tilde{\Omega}_{A\setminus S(m),A}(I-\Lambda)^{-1}_{A,P(m)}\alpha
\\
&=&\tilde{\Omega}_{A\setminus
S(m),[m]}(I-\Lambda)^{-1}_{[m],P(m)}\alpha\\
&&{}-
\tilde{\Omega}_{A\setminus S(m),[m]\setminus A}
(I-\Lambda)^{-1}_{[m]\setminus A,P(m)}\alpha.
\end{eqnarray*}
By hypothesis, the first term in the last line is zero. By Lemma~\ref
{lem:sublemma1}
below, $(I-\Lambda)^{-1}_{[m]\setminus A,P(m)}\alpha=0$, and so the
second term in the last line is zero as well. Therefore,
\[
\tilde{\Omega}_{A\setminus
S(m),A}(I-\tilde\Lambda)^{-1}_{A,\tilde{P}(m)}\alpha_{\tilde{P}(m)}=0.
\]

It remains to be shown that $\alpha_{\tilde{P}(m)}\neq0$. Suppose
instead that $\alpha_{\tilde{P}(m)}= 0$. Then, using
Lemma~\ref{lem:sublemma1}, we obtain that
\begin{eqnarray*}
0&=&(I-\Lambda)^{-1}_{[m]\setminus A,P(m)}\alpha\\
&=&(I-\Lambda)^{-1}_{[m]\setminus
A,\tilde{P}(m)}\alpha_{\tilde{P}(m)}+(I-\Lambda
)^{-1}_{[m]\setminus
A,P(m)\setminus\tilde{P}(m)}\alpha_{P(m)\setminus
\tilde{P}(m)}\\
&=&0+(I-\Lambda)^{-1}_{[m]\setminus A,P(m)\setminus
\tilde{P}(m)}\alpha_{P(m)\setminus\tilde{P}(m)}.
\end{eqnarray*}
However, $P(m)\setminus\tilde{P}(m)\subseteq[m]\setminus A$ and thus\vspace*{-1pt}
$(I-\Lambda)^{-1}_{[m]\setminus A,P(m)\setminus\tilde{P}(m)}$ is a
submatrix of $(I-\Lambda)^{-1}_{[m]\setminus A,[m]\setminus A}$, which
is a full rank matrix as it is upper triangular with ones on the
diagonal. Therefore, $(I-\Lambda)^{-1}_{[m]\setminus A,P(m)\setminus
\tilde{P}(m)}$ is\vspace*{1pt} full rank, and so $\alpha_{P(m)\setminus
\tilde{P}(m)}=0$. It follows that $\alpha=0$, which is a
contradiction. We conclude that $\alpha_{\tilde{P}(m)}\neq0$ and, by
Lemma~\ref{lem:overall-ident}, that $\phi_{G_W}$ is not injective.
\end{pf}

\begin{proposition}
\label{prop:sufficiency}
Let $G=(V,D,B)$ be an acyclic mixed graph with topologically ordered
vertex set $V=[m+1]$, such that the parametrization $\phi_G$ is not
injective. If either the directed part of $G$ does not contain an
arborescence converging to $m+1$, or the bidirected part of $G$ is not
connected, then there is some proper induced subgraph $G_W$ of $G$ for
which the parametrization $\phi_{G_W}$ is not injective.
\end{proposition}

\begin{pf}
From Lemma~\ref{lem:overall-ident}, for some $i\leq m$, $\Lambda\in
\mathbb{R}^D$  and $\Omega\in\mathit{PD}(B)$,
%
%e5.1 ###
\begin{equation}
\label{eq:suff-rank}
\rank \bigl( \Omega_{[i]\setminus S(i),
[i]}(I-\Lambda)^{-1}_{[i],P(i)}  \bigr)<|P(i)| .
\end{equation}
Suppose $i<m$. Take $W=[i+1]$, and denote the induced subgraph\vspace*{-1pt} as
$G_W=(W,\tilde D,\tilde B)$. It holds trivially that
$\tilde{\Lambda}:=\Lambda_{[i+1],[i+1]}\in\mathbb{R}^{\tilde D}$ and
$\tilde{\Omega}:=\Omega_{[i+1],[i+1]}\in\mathit{PD}(\tilde B)$, and
furthermore
$(I-\tilde{\Lambda})^{-1}=(I-\Lambda)^{-1}_{[i+1],[i+1]}$. It
is then clear that, by Lemma~\ref{lem:overall-ident}, $\phi_{G_W}$ is
not injective.

Next suppose instead that (\ref{eq:suff-rank}) is true for $i=m$. If the
directed part of $G$ does not contain an arborescence converging to
$m+1$, then apply Proposition~\ref{prop:directed_sufficiency} to produce
a proper induced subgraph $G_W$ with $\phi_{G_W}$ noninjective. If
instead the bidirected part of $G$ is not connected, then apply
Proposition~\ref{prop:bidirected_sufficiency} to produce a proper induced
subgraph $G_W$ with $\phi_{G_W}$ noninjective.

In all cases, we have constructed a subset $W\subsetneq V$ with
$\phi_{G_W}$ not injective.
\end{pf}

%s6 ###
\section{\texorpdfstring{Proofs of lemmas in Section~\protect\ref{sec:sufficiency}}%
{Proofs of lemmas in Section 5}}
\label{sec:lemma_proofs}

%% Here we prove lemmas used in Section~\ref{sec:sufficiency}.

\begin{lemma}
\label{lem:sublemma1}
Let $G$, $\Lambda$, $\Omega$, $\alpha$, and $A$ be as in the
statement of Proposition~\ref{prop:bidirected_sufficiency}.
Then $(I-\Lambda)^{-1}_{[m]\setminus A,P(m)}\alpha=0$.
\end{lemma}

\begin{pf}
If $i\in[m]\setminus A$ and $j\in A$, then, by definition of $A$,
it holds that \mbox{$\Omega_{i,j}=0$}. Therefore, $\Omega_{[m]\setminus
A,A}=0$ and we obtain that
\[
\Omega_{[m]\setminus A,[m]\setminus A}
(I-\Lambda)^{-1}_{[m]\setminus A,P(m)}\alpha=
\Omega_{[m]\setminus A,[m]}(I-\Lambda)^{-1}_{[m],P(m)}\alpha=0.
\]
For the last equality, observe that $[m]\setminus A\subset
[m]\setminus S(i)$ since $S(i)\subset A$. Since
$\Omega_{[m]\setminus A,[m]\setminus A}$ is positive definite,
the claim follows.
\end{pf}

For a directed path $\pi$ in the graph $G$, we write $\pi\not\subset G_A$
to indicate that not all the nodes of $\pi$ lie in $A$. Also, by
convention, $\mathcal{P}(j,j)$ is a singleton set containing the trivial
path at $j$; in this case $\pi$ has no edges and we define $\prod
_{a\to
b\in\pi}\lambda_{ab}=1$.

\begin{lemma}
\label{lem:sublemma2}
Let $G$, $\Lambda$, $\Omega$, $\alpha$, and $A$ be as in the
statement of
Proposition~\ref{prop:bidirected_sufficiency}. Then for every $i\leq
m$,
\[
\sum_{k\in P (m )}\alpha_k \biggl(\sum_{\pi\in
\mathcal{P} (i,k ),\pi\not\subset G_A}
\prod_{a\to b\in\pi}\lambda_{ab} \biggr)=0.
\]
\end{lemma}

\begin{pf}
First, we prove the claim for $i\notin A$. Working from
Lemma~\ref{lem:sublemma1}, we have that\vspace{-4pt}
%
%e6.1 ###
\begin{eqnarray}
\label{eq:sublemma2-1}
0&=&(I-\Lambda)^{-1}_{i,P(m)}\alpha=\sum_{k\in
P(m)}(I-\Lambda)^{-1}_{ik}\alpha_k\nonumber
\\[-10pt]
\\[-10pt]
&=&\sum_{k\in P(m)}\alpha_k \biggl(\sum_{\pi\in\mathcal{P}(i,k)}
\prod_{a\to b\in\pi}\lambda_{ab} \biggr).
\nonumber
\end{eqnarray}
Since $i\notin A$, any path $\pi\in\mathcal{P}(i,k)$ for any $k$
necessarily satisfies $\pi\not\subset G_A$. Hence, we can rewrite
(\ref{eq:sublemma2-1}) as\vspace{-2pt}
\[
\sum_{k\in P(m)}\alpha_k \biggl(
\sum_{\pi\in\mathcal{P}(i,k),\pi\not\subset G_A}
\prod_{a\to b\in\pi}\lambda_{ab} \biggr) =0.
\]

Next, we address the case $i\in A$. Inducting on $i$ in decreasing order,
we may assume that the claim holds for all $j\in\{i+1,i+2,\ldots ,m\}$.
[As a base case, we can set $i=m$ because, by the assumed topological
order, $\mathcal{P}(m,k)=\varnothing$ for all nodes $k<m$.] The quantity
claimed to be vanishing is
%
%e6.2 ###
\begin{eqnarray}
\label{eq:sublemma2-2}
&&\sum_{k\in
P(m)}\alpha_k \biggl(\sum_{\pi\in\mathcal{P}(i,k),\pi\not\subset
G_A} \prod_{a\to b\in\pi}\lambda_{ab} \biggr) \nonumber
\\[-10pt]
\\[-10pt]
&& \qquad =\sum_{k\in P(m)}\alpha_k \biggl[\sum_{j\dvtx i\to
j} \biggl(\sum_{\pi'\in\mathcal{P}(j,k),\pi'\not\subset G_A}
\lambda_{ij}\prod_{a\to b\in\pi'}\lambda_{ab} \biggr) \biggr].
\nonumber
\end{eqnarray}
This last equality is obtained by splitting any path $\pi=i\to v_1\to
\cdots \to v_n\to k$ into $i\to j:=v_1$ and $\pi'=j\to v_2\to\cdots \to
v_n\to k$. (Note that the path of length zero at $i$ is not in the sum,
since this path would not satisfy $\pi\not\subset G_A$.)
Since we assume $i\in A$, it holds that $\pi\not\subset G_A$
if and only if $\pi'\not\subset G_A$. Interchanging the order of the
summations in (\ref{eq:sublemma2-2}), we obtain that
\begin{eqnarray*}
&&\sum_{k\in
P(m)}\alpha_k \biggl(\sum_{\pi\in\mathcal{P}(i,k),\pi\not\subset
G_A} \prod_{a\to b\in\pi}\lambda_{ab} \biggr) \\[-2pt]
&& \qquad =\sum_{j\dvtx i\to j} \biggl[\sum_{k\in
P(m)}\alpha_k \biggl(\sum_{\pi'\in\mathcal{P}(j,k),\pi'\not
\subset
G_A} \lambda_{ij}\prod_{a\to b\in
\pi'}\lambda_{ab} \biggr) \biggr]\\[-2pt]
&& \qquad =\sum_{j\dvtx i\to
j}\lambda_{ij} \biggl[\sum_{k\in P(m)}\alpha_k
 \biggl(\sum_{\pi'\in\mathcal{P}(j,k),\pi'\not\subset G_A}
\prod_{a\to b\in\pi'}\lambda_{ab} \biggr) \biggr].
\end{eqnarray*}
Working with a topologically ordered set of nodes, the presence of
an edge $i\to j$ implies $i<j$. The inductive hypothesis thus
yields that
\[
\sum_{k\in
P(m)}\alpha_k  \biggl(\sum_{\pi\in\mathcal{P}(i,k),\pi\not\subset
G_A} \prod_{a\to b\in\pi}\lambda_{ab} \biggr)=\sum_{j\dvtx i\to
j}\lambda_{ij}\cdot0=0,
\]
which completes the inductive step and the proof of the lemma.
\end{pf}

\begin{lemma}
\label{lem:sublemma3}
Let $G$, $\Lambda$, $\Omega$, $\alpha$ and $A$ be as in the
statement of Proposition~\ref{prop:bidirected_sufficiency}. Then
for all $i\in A$,
\[
(I-\tilde{\Lambda})^{-1}_{i,\tilde{P}(m)}\alpha_{\tilde{P}(m)}=
(I-\Lambda)^{-1}_{i,P(m)}\alpha.
\]
\end{lemma}

\begin{pf}
The right-hand side of the above equation can be rewritten as
\begin{eqnarray*}
 (I-\Lambda)^{-1}_{i,P(m)}\alpha&=&\sum_{k\in
P(m)}(I-\Lambda)^{-1}_{ik}\alpha_k=\sum_{k\in
P(m)}\alpha_k \biggl(\sum_{\pi\in\mathcal{P}(i,k)} \prod_{a\to
b\in\pi} \lambda_{ab} \biggr)\\
  &=&\sum_{k\in
P(m)}\alpha_k \biggl(\sum_{\pi\in\mathcal{P}(i,k),\pi\subset
G_A}\prod_{a\to b\in\pi} \lambda_{ab} \biggr)\\
&&{} +\sum_{k\in
P(m)}\alpha_k \biggl(\sum_{\pi\in\mathcal{P}(i,k),\pi\not\subset
G_A}\prod_{a\to b\in\pi} \lambda_{ab} \biggr).
\end{eqnarray*}
Consider the two sums in the last line above.
By Lemma~\ref{lem:sublemma2}, the second sum is equal to zero. Note
also that if $k\in P(m)\setminus A$, then there is no path
$\pi\in\mathcal{P}(i,k)$ with $\pi\subset G_A$. Therefore, the first
sum can be indexed over $k\in\tilde{P}(m)$. We thus obtain that, as
claimed,
\begin{eqnarray*}
 (I-\Lambda)^{-1}_{i,P(m)}\alpha&=&\sum_{k\in
\tilde{P}(m)}\alpha_k \biggl(\sum_{\pi\in\mathcal{P}(i,k),\pi
\subset
G_A}\prod_{a\to b\in\pi} \lambda_{ab} \biggr)\\
%%=\sum_{k\in P(m)\cap A}
  &=&\sum_{k\in\tilde{P}(m)}
\alpha_k(I-\tilde{\Lambda})^{-1}_{ik}
=(I-\tilde{\Lambda})^{-1}_{i,\tilde{P}(m)}\alpha_{\tilde{P}(m)}.
\end{eqnarray*}
\upqed
\end{pf}

%s7 ###
\section{Cyclic models}
\label{sec:cyclic-models}

In this section, we prove Theorem~\ref{thm:acyclic} from the\break\hyperref[sec:intro]{Introduction},
which states that only acyclic mixed graphs may yield globally identifiable
models. By Lemma~\ref{lem:simple}, the theorem holds if we can show that
the parametrization $\phi_G$ is not injective when $G$ is a simple directed
cycle, that is, when $G$ is isomorphic to the cycle
%
%e7.1 ###
\begin{equation}
\label{eq:cycle-ordered}
1\to2\to\cdots \to m\to1
\end{equation}
for some $m\ge3$. This noninjectivity is shown in the next lemma.
Recall the definition of a fiber in (\ref{eq:fiber}).

\begin{lemma}
\label{lem:cycles}
Let $G=(V,D,B)$ be a simple directed cycle on $m\ge3$ nodes,
$\Lambda\in\mathbb{R}^D_{\mathrm{reg}}$ and $\Omega\in\mathit
{PD}(B)$. Then
the cardinality of the fiber $\mathcal{F}(\Lambda,\Omega)$ is at
most two
and is equal to two for generic choices of $\Lambda$ and $\Omega$.
\end{lemma}

In order to prepare the proof of Lemma~\ref{lem:cycles}, note that for
directed graphs the set $\mathit{PD}(B)=\mathit{PD}(\varnothing)$ contains
exactly the diagonal matrices with positive diagonal entries. This set
being invariant under matrix inversion, it is convenient to consider
the polynomial map
\[
\kappa_G\dvtx  (\Lambda,\Delta) \mapsto(I-\Lambda)\Delta(I-\Lambda)^T
\]
that parametrizes the inverse of the covariance matrix of the distributions
in the structural equation model. Since
$\kappa_G(\Lambda,\Delta)=\phi_G(\Lambda,\Delta^{-1})^{-1}$ for
$\Lambda\in\mathbb{R}^D_{\mathrm{reg}}$ and $\Delta\in\mathit
{PD}(\varnothing)$,
the fibers of $\kappa_G$ and $\phi_G$ are in bijection with each other.

\begin{pf*}{Proof of Lemma~\ref{lem:cycles}}
Without loss of generality, assume $G$ to be the graph with the edges in
(\ref{eq:cycle-ordered}). For shorter notation, we let
$\lambda_i=\Lambda_{i,i+1}$, the parameter on the edge $i\to i+1$.
Throughout, indices are read cyclically with $m+i:=i$ for $i\ge1$. The
matrix $(I-\Lambda)$ is invertible if and only if $\prod_{i=1}^m
\lambda_i\neq1$. Let $\delta_i=\Delta_{ii}$, the inverse of the
positive variance parameter associated with node $i$. Treating
$\kappa_G$ as a function of a pair of vectors
$(\lambda,\delta)\in\mathbb{R}^m\times\mathbb{R}_+^m$, we obtain that
$\kappa_G(\lambda,\delta)$ is equal to
\[
\pmatrix{
\delta_1+\delta_2\lambda_1^2&-\delta_2\lambda_1&0&\cdots &0&-\delta
_1\lambda_m\cr
-\delta_2\lambda_1&\delta_2+\delta_3\lambda_2^2&-\delta_3\lambda
_2&\cdots &0&0\cr
0&-\delta_3\lambda_2&\delta_3+\delta_4\lambda_3^2&\cdots &0&0\cr
\cdots &\cdots &\cdots &\cdots &\cdots &\cdots \cr
0&0&0&\cdots &\delta_{m-1}+\delta_m\lambda_{m-1}^2&-\delta_m\lambda
_{m-1}\cr
-\delta_1\lambda_m&0&0&\cdots &-\delta_m\lambda_{m-1}&\delta
_m+\delta_1\lambda_m^2
}
.
\]

Fix a pair $(\lambda^0,\delta^0)\in\mathbb{R}^m\times\mathbb
{R}_+^m$ with
$\prod_{i=1}^m
\lambda_i^0\neq1$. We
wish to describe the fiber
%
%e7.2 ###
\begin{equation}
\label{eq:fiber-inv}
%%\mathcal{H}(\lambda^0,\delta^0) =
 \{
(\lambda,\delta)\in\mathbb{R}^m\times\mathbb{R}_+^m \dvtx
\kappa_G(\lambda,\delta)=\kappa_G(\lambda^0,\delta^0)
 \}.
\end{equation}
Let $K^0:=\kappa_G(\lambda^0,\delta^0)$. The equation
$\kappa_G(\lambda,\delta)=K^0$ determining membership in the fiber
amounts to the system of the $2m$ polynomial equations
%% we enumerate as
%% (\ref{eqn:i.1}) and (\ref{eqn:i.2}) for $i=1,2,\dots,m$:
\renewcommand{\theequation}{7.3a.$i$}
\begin{eqnarray}
\label{eqn:i.1}
 \delta_i+\delta_{i+1}\lambda
_i^2=K^0_{i,i} ,\\[-37pt]\nonumber
\end{eqnarray}
 \renewcommand{\theequation}{7.3b.$i$}
\begin{equation}
\label{eqn:i.2}
 \hspace*{15pt}-\delta_{i+1}\lambda_i=K^0_{i,i+1}
\end{equation}
for $i=1,\ldots ,m$. We split the problem into two cases, for which the
algebraic degree of the equation system given by (\ref{eqn:i.1}) and
(\ref{eqn:i.2}) differs.

\textit{Case (i)}: Suppose $\lambda^0_i=0$ for some $i$. Without loss of
generality, $\lambda^0_1=0$ such that $K^0_{12}=0$ and
$K^0_{11}=\delta^0_1$. As a consequence, the two equations
(\ref{eqn:i.1}) and (\ref{eqn:i.2}) for $i=1$ reduce to
$\delta_1=\delta^0_1$ and $\lambda_1=0=\lambda^0_1$. This provides the
basis for solving the remaining equations recursively in the order
$i=m,\ldots ,2$. Each time the equation pair reduces to the linear
equations $\delta_i=\delta^0_i$ and $\lambda_i=\lambda^0_i$, and the
fiber in (\ref{eq:fiber-inv}) is seen to be the singleton
$\{(\lambda^0,\delta^0)\}$. Note that the problem has become the same as
parameter identification in the model based on the acyclic graph obtained
by removing the edge $1\to2$ from~$G$. Note further that the equation
system is of degree one in this case.

\textit{Case (ii)}: Assume now that $\lambda^0_i\neq0$ for all $i$. We
claim that the fiber in (\ref{eq:fiber-inv}) then also contains the pair
$(\lambda^1,\delta^1)$ that has coordinates
\begin{eqnarray*}
\delta^1_i&=&\delta^0_i+
\frac{ (\prod_{j=1}^m
\delta^0_j )
 [\prod_{j=1}^m
 ( (\lambda^0_j )^2-1 ) ]
}{\det (K^0_{-i} )}, \\
\lambda^1_i&=&-\frac{K^0_{i,i+1}}{\delta_{i+1}}
\end{eqnarray*}
for $i=1,2,\ldots ,m$. Here $K^0_{-i}$ is the matrix obtained from $K^0$
by removing the $i$th row and column. Note that
$(\delta^1,\lambda^1)\not=(\delta^0,\lambda^0)$ if and only if
$\prod_{j=1}^m \lambda^0_j\neq-1$; recall that the product is
assumed to
be different from $1$ to ensure that $I-\Lambda$ is invertible. It is
not very difficult to check that $(\delta^1,\lambda^1)$ is indeed in the
fiber; the $m$ equations in (\ref{eqn:i.2}) are satisfied trivially, and
the $m$ equations in (\ref{eqn:i.1}) can be checked by plug-in. For this
an explicit expression of $\det (K^0_{-i} )$ in terms of
$(\lambda^0,\delta^0)$ is needed. Using the Cauchy--Binet formula, one
can show that
\[
\det (K^0_{-i} ) =
 \Biggl(\prod_{j=1}^m \delta^0_j \Biggr)
 \Biggl( \frac{1}{\delta^0_i} + \sum_{j=1}^{i-1}
\frac{1}{\delta^0_j}\prod_{k=j}^{i-1}
 (\lambda^0_k )^2
+ \sum_{j=i+1}^{m} \frac{1}{\delta^0_j}\prod_{k=j}^{m+i-1}
 (\lambda^0_k )^2 \Biggr).
\]

We furthermore claim that the fiber contains no points other than
$(\lambda^0,\delta^0)$ and $(\lambda^1,\delta^1)$. We outline the proof
of this claim, again leaving out some of the details.

Solve for $\lambda_1$ in equation (\ref{eqn:i.2}) for $i=1$ and plug the
resulting expression in $\delta_2$ into the equation (\ref{eqn:i.1}) for
$i=1$. This equation can be solved for $\delta_2$ to give an expression
in $\delta_1$. Continue on in this fashion for the indices $i=2,\ldots ,m$
always obtaining an expression in $\delta_1$ after solving
(\ref{eqn:i.1}). Let $[j\dvtx k]:=\{j,\ldots ,k\}$ for integers $j<k$. We
find that, after the $i$th step,
\[
\delta_i = (K^0_{i-1,i})^2 \cdot
\frac{\det (K^0_{[1\dvtx i-2],[1\dvtx i-2]} )-
\det (K^0_{[2\dvtx i-2],[2\dvtx i-2]} )\delta_1}{
\det (K^0_{[1\dvtx i-1],[1\dvtx i-1]} )-
\det (K^0_{[2\dvtx i-1],[2\dvtx i-1]} )\delta_1},
\]
where we define
$\det(K^0_{[1\dvtx 0]})=\det(K^0_{[2\dvtx 1]})=1$ and
$\det(K^0_{[2\dvtx 0]})=0$. The last step of this procedure,
namely, plugging the expression for $\delta_m$ into the equation
(\ref{eqn:i.1}) for $i=m$ produces a rational equation in the single
variable $\delta_1$. Clearing denominators we obtain a quadratic
equation in $\delta_1$ whose leading coefficient for $\delta_1^2$
simplifies to $\det(K^0_{-1})$ and thus is nonzero. Therefore, the
polynomial equation system in (\ref{eqn:i.1})--(\ref{eqn:i.2}) has degree
two and the fiber in (\ref{eq:fiber-inv}) contains precisely
$(\lambda^0,\delta^0)$ and $(\lambda^1,\delta^1)$. Note that the fiber
has cardinality one (with a point of multiplicity two) if $\prod_{j=1}^m
\lambda^0_j= -1$.
\end{pf*}

%s8 ###
\section{Conclusion}
\label{sec:conclusion}

Our Theorems~\ref{thm:acyclic} and \ref{thm:main} fully characterize the
mixed graphs for which the associated linear structural equation model is
globally identifiable. Globally identifiable models have smooth manifolds
as parameter spaces, which implies in particular that maximum likelihood
estimators are asymptotically normal for all choices of a true distribution
in the model. Similarly, likelihood ratio statistics for testing two
nested globally identifiable models are asymptotically chi-square.
Example~\ref{ex:simulation} demonstrates that these properties may
fail in
models that are only generically identifiable. The resulting inferential
issues are also not so easily overcome using bootstrap methods; compare
\cite{andrews2010}. Nevertheless, generically identifiable models
appear in
various applications, and characterizing the mixed graphs that yield
generically identifiable linear structural equation models remains an
important open problem.

\section*{Acknowledgments}
We are grateful to two referees and an associate editor who provided very
helpful comments on the original version of this paper.

%suskaldyti doi

\printaddresses

\end{document}